\documentclass[review,hidelinks,onefignum,onetabnum]{siamart220329}
\usepackage{amssymb,amsmath} 
\usepackage{subeqnarray,comment} 
\usepackage{graphicx}
\usepackage{placeins}
\usepackage{float}
\usepackage{array}
\usepackage{subfigure}
\usepackage{indentfirst}
\usepackage{multicol}
\usepackage{multirow}
\newcommand{\Reals}[1]{{\rm I\! R}^{#1}}
\newtheorem{Theorem}{Theorem}
\nolinenumbers

\usepackage{lipsum}
\usepackage{amsfonts}
\usepackage{graphicx}
\usepackage{epstopdf}
\usepackage{algorithmic}
\ifpdf
  \DeclareGraphicsExtensions{.eps,.pdf,.png,.jpg}
\else
  \DeclareGraphicsExtensions{.eps}
\fi

\newsiamremark{remark}{Remark}
\newsiamremark{hypothesis}{Hypothesis}
\crefname{hypothesis}{Hypothesis}{Hypotheses}
\newsiamthm{claim}{Claim}

\headers{An optimal $H^1$ norm error estimate for Finite Neuron Method}{J. Jia, Y. Lee and R. Shan}

\title{An unconstrained formulation of some constrained partial differential equations and its application to finite neuron methods\thanks{Submitted to the editors DATE.
\funding{The first author was funded in part by NSFC 22341302 and the second and corresponding author was funded in part by NSF-DMS 2208499}}}

\author{Jiwei Jia\thanks{School of Mathematics, Jilin University, Changchun 130012, Jilin,  China 
 (\email{jiajiwei@jlu.edu.cn}, \email{shanrt22@mails.jlu.edu.cn}).}
\and Young Ju Lee\thanks{Department of Mathematics, Texas State University, San Marcos, Texas, USA 
  (\email{yjlee@txstate.edu}).}
\and Ruitong Shan\footnotemark[2]
}

\usepackage{amsopn}

\ifpdf
\hypersetup{
  pdftitle={An unconstrained formulation of some constrained partial differential equations and its application to finite neuron methods},
  pdfauthor={Jiwei Jia, Young Ju Lee and Ruitong Shan}
}
\fi

\begin{document}

\maketitle

\begin{abstract}
In this paper, we present a new framework how a PDE with constraints can be formulated into a sequence of PDEs with no constraints, whose solutions are convergent to the solution of the PDE with constraints. This framework is then used to build a novel finite neuron method to solve the 2nd order elliptic equations with the Dirichlet boundary condition. Our algorithm is the first algorithm, proven to lead to shallow neural network solutions with an optimal $H^1$ norm error. We show that a widely used penalized PDE, which imposes the Dirichlet boundary condition weakly can be interpreted as the first element of the sequence of PDEs within our framework. Furthermore, numerically, we show that it may not lead to the solution with the optimal $H^1$ norm error bound in general. On the other hand, we theoretically demonstrate that the second and later elements of a sequence of PDEs can lead to an adequate solution with the optimal $H^1$ norm error bound. A number of sample tests are performed to confirm the effectiveness of the proposed algorithm and the relevant theory. 
\end{abstract}



\section{Introduction}

There are many PDEs with constraints in real-life applications. Among others, Darcy's law model, Stokes and Navier-Stokes equations contain the continuity equation as constraints. Typically, such PDEs with constraints are solved using the framework of saddle point formulations, where the Lagrange multipliers are introduced to satisfy the constraints. Generally, the first step to solve the saddle point problem is to discretize the equations for both the primary variable and the Lagrange multiplier. The most difficult issue in this approach lies at ensuring the well-posedness of the discrete system, namely, the Babuska-Brezzi or discrete inf-sup condition \cite{babuvska1973finite,fortin1991mixed,xu2003some}. In general, the establishment of discrete inf-sup condition is nontrivial. 

In this paper, we present a novel approach to handle a class of PDEs with constraints in a way that the solution can be approximated via a solution to the unconstrained PDE. Such an approach can be shown to be valid if the original PDE with constraints is well-posed. The main idea lies at applying the augmented Lagrangian Uzawa iteration \cite{lee2007robust}
to solve the constrained PDEs in purely continuous level. The augmented Lagrangian Uzawa iteration introduces a sequence of PDEs to be solved, which are unconstrained ones and the corresponding solution can be made to be arbitrarily close within a few iterations, to the solution to that of the original PDEs of interest by choosing an appropriate penalty parameter. We remark that in literature, some penalized PDEs have been used to handle the PDE with constraints, which is traced back to Lion's formulation \cite{lions2010problemes} for a special case of the 2nd order elliptic PDE with the Dirichlet boundary condition. We demonstrate that it is the special case of the first element of the sequence of PDEs obtained from the augmented Lagrangian Uzawa iteration within our framework. 

The main motivation behind such an unconstrained formulation of the PDEs with constraints is at attempting to resolve an open issue in the solution to the 2nd order elliptic equation with the Dirichlet boundary condition using the shallow neural networks \cite{xu2020finite}. When the neural networks are used as an approximation function space to the PDEs with Dirichlet boundary condition, it is a common practice to impose the Dirichlet boundary condition as a penalty due to the global nature of the neural network functions. Namely, it is unlike finite element methods, where the local nature of the finite element basis allows the strong imposition of the dirichlet boundary condition. The use of the a penalty formulation for enforcing the Dirichlet boundary condition can be traced back to the late 1960s \cite{lions2010problemes}, where the elliptic problems with very rough Dirichlet boundary data was considered, namely, 
\begin{equation}
-\Delta u = f \mbox{ in } \Omega \quad \mbox{ and } \quad u = g \quad \mbox{ on } \partial \Omega
\end{equation}
and the regularization of the problem was considered by replacing the Dirichlet boundary conditions with the Robin boundary condition (see also \cite{arnold2002unified}), i.e., 
\begin{equation}
-\Delta u_\delta = f \quad \mbox{ and } \quad u_\delta + \delta^{-1} \frac{\partial u_\delta}{\partial n} = g, 
\end{equation}
where $\delta^{-1}$ is the penalty parameter. Furthermore, it was proven that when $\delta$ approaches $0$, the solution $u_\delta$ approaches $u$.


Such a penalty method has a particular advantage since the inf-sup condition is not needed in the discrete setting and the shallow neural network solution can be trained using the orthogonal greedy algorithm \cite{siegel2021greedy} since the penalized PDE is an elliptic problem. This approach is actually used to obtain neural network solutions to the PDEs (see \cite{yu2018deep,xu2020finite,raissi2017physics,cai2020deep} and references cited therein). In particular, the finite neuron method introduced in \cite{xu2020finite} has also used the Lion's formulation \cite{lions2010problemes} to solve the 2nd order elliptic PDEs with the Dirichlet boundary conditions. On the other hand, the optimal error estimate in $H^1$ norm was not established. Namely, it was shown that the guaranteed error estimate in this framework can be given as follows: 
\begin{equation}\label{eq:estimate}
\|u - u_n\|_{H^1} \lesssim n^{-\frac{1}{4} - \frac{2(k-1)+1}{4d}}, 
\end{equation}
where ${\rm ReLu}^k$ activation function is used for the shallow neural network functions and $d$ is the dimension. Furthermore, the error estimate \eqref{eq:estimate} was established under the condition that the penalty parameter $\delta^{-1}$ is as large as $n^{-\frac{1}{2}-\frac{2(k-1)+1}{2d}}$. 
The main result in this paper is a new algorithm, which can lead to the optimal error estimate \cite{siegel2022sharp}: 
\begin{equation}\label{main:est} 
\|u - u_n\|_{H^1} \lesssim n^{-\frac{1}{2} - \frac{2(k-1)+1}{2d}}. 
\end{equation}
More importantly, we show that the error estimate \eqref{main:est} for the newly proposed algorithm can be obtained under the condition that the penalty paramter $\delta^{-1}$ is just as large as $n^{-\frac{1}{3}-\frac{2(k-1)+1}{3d}}$.  

The new algorithmic development lies at first reformulating the PDEs with the Dirichlet boundary condition as PDEs with constraints. The latter is then formulated as the saddle point problem. We then consider the augmented Lagrangian Uzawa iteration to solve the saddle point problem. Basically, this framework leads to a sequence of unconstrained elliptic equations to be solved. We shall show that while the first element of the sequence of unconstrained PDEs may not be able to lead to the optimal solution. The second or later elements can lead to the optimal neural network solutions if the penalty parameter $\delta^{-1}$ is chosen appropriately (see \S \ref{num:ex} for numerical experiments). In particular, we demonstrate the theoretically  optimal convergence in $H^1$ norm for the second element of the sequence of unconstrained PDEs. We remark that even if the optimal accuracy can be guaranteed, the disadvange exists, which is still to solve the first element of the sequence to obtain the solution to the second element of the sequence. On the other hand, unlike the penalized PDEs pioneered at \cite{xu2020finite}, the optimal accuracy can be attained with relatively small penalty parameter $\delta^{-1}$. It is obvious to see that the larger the penalty parameter is, the more difficult it is to apply the orthogonal greedy algorithm. Again, the proposed algorithm requires the solution to the decoupled simple systems, which can be handled by the orthogonal Greedy algorithms to train the shallow neural network solution, thus the construction of the approximate solutions does not need the discrete inf-sup condition. For the technical tool to obtain the $H^1$ error estimate, we employ the trace theorem, which results in $L^2$ and $H^1$ norm error in a way that the penalty parameter is drawn to the $L^2$ norm error. We then apply the duality argument to estimate the resulting $L^2$ norm error, following Shatz and Wang \cite{schatz1996some}. This is shown to handle the penalty parameter in an appropriate manner. The similar idea has been successfully used in authors' recent work to approximate the indefinite problem via the shallow neural network functions \cite{hongjialeeli2024}. Basically, the duality argument lessens the difficulty due to the penality parameter in obtaining the optimal error estimate.  


Throughout the paper, we shall use the standard notation. Following \cite{xu1992iterative}, we use $A \lesssim B$ and $A \gtrsim B$ if there is a generic constant $c$ such that $A \leq c B$ and $A \geq c B$, respectively. As an approximation space of the solution to PDEs, a stable shallow neural network function space is recommended to be used \cite{siegel2021greedy} and a stable shallow neural network function space is defined by the following: for a fixed positive reals $B$ and $M$, and an integer $n$,  
\begin{equation}
\Sigma_{d}^{n}(\sigma) = \left \{ \sum_{i=1}^n a_i \sigma_k(\omega_i \cdot x + b_i) : \omega_i \in \mathbb{S}^{d-1}, \, |b_i| \leq B, \, \sum_{i=1}^n |a_i| \leq M \right \},  
\end{equation}
where $\mathbb{S}^{d-1}$ is the unit sphere in $\Reals{d}$ and $\sigma_k$ is the ReLu$^k$ function. It was shown that the space $\Sigma_d^n(\sigma)$ has the following approximation property \cite{siegel2022sharp,siegel2021greedy}:
\begin{equation}
\inf_{u_n \in \Sigma^{n}_d(\sigma)} \|u - u_n\|_{H^m} \lesssim \left ( n^{-\frac{1}{2} - \frac{2(k-m)+1}{2d}}\right ),   
\end{equation}
where $H^m$ denotes the standard Sobolev space of $m$th order. Note that in this paper, we shall focus our concern to the case when $m = 1$. Throughout the paper, we let $\Omega \in \Reals{d}$ (with $d = 1,2,3$) be the smooth domain with the boundary $\Gamma$. The Sobolev spaces $H^s(K)$ for $K \subset \Omega$ or $K \subset \Gamma$ and $s \geq 0$ are defined in the standard way \cite{adams2003sobolev}. The standard $L^2$ space will be used as well. The Sobolev norms and seminorms will be denoted by $\|\cdot\|_{s,K}$ and $|\cdot|_{s,K}$ with the subscript $K$ dropped when $K = \Omega$, respectively. We recall that the following trace inequality holds for $v \in H^s(\Omega)$ with $s > \frac{1}{2}$, 
\begin{equation}
\|v\|_{s - \frac{1}{2}} \leq C \|v\|_s.
\end{equation}
We shall also use the space $H^{-\frac{1}{2}}(\Gamma)$, i.e., the dual space of $H^{\frac{1}{2}}(\Gamma)$ with the following norm:
\begin{equation}
\|\mu\|_{-\frac{1}{2}} = \sup_{z \in H^{\frac{1}{2}}(\Gamma)} \frac{\langle \mu, z\rangle}{\|z\|_{\frac{1}{2}}}. 
\end{equation}
Note that $\langle \cdot,\cdot \rangle$ is the dual pairing. For functions $v \in H^1(\Omega)$ with $\Delta v \in L^2(\Omega)$, it holds that 
\begin{equation}
\| \nabla v \cdot n \|_{-\frac{1}{2}} \lesssim \|v\|_1 + \|\Delta v\|_0. 
\end{equation}
For a given $g \in H^{\frac{1}{2}}(\Gamma)$, the symbol $H^1_g(\Omega)$ denotes  the subset of $H^1(\Omega)$ defined as follows:
\begin{equation}
H^1_g(\Omega) = \left \{ u \in H^1(\Omega) : {\rm tr}(u) = g \, \mbox{ on } \Gamma \right \},   
\end{equation}
where ${\rm tr}$ is the standard trace operator defined on $H^1(\Omega)$ (see \cite{adams2003sobolev} and \cite{girault2012finite}). Lastly, we shall denote the dual of $H^1(\Omega)$ by $H^{-1}(\Omega)$. Lastly, for any given operator, $O$, ${\rm ker}(O)$ denotes the null space of the operator $O$. 

\section{Governing Equations}
Given a smooth domain $\Omega$ and its boundary $\Gamma$, we shall consider to solve the elliptic equation with the Dirichlet boundary condition given as follows: 
\begin{subeqnarray}\label{eq:main}
-\Delta u + a_0 u = f_0, \quad \mbox{ in } \Omega, \\ 
u = g, \quad \mbox{ on } \Gamma,   
\end{subeqnarray}
where $a_0 \geq 0$, $f_0 \in H^{-1}(\Omega)$ and $g \in H^{\frac{1}{2}}(\Gamma)$. We remark that our analysis allows $a_0 = 0$. This equation can be formulated as an optimization problem given as follows:
\begin{equation}\label{eq:orig}
u = {\rm arg} \left \{ \min_{v \in H_g^1(\Omega)} \frac{1}{2} \int_\Omega |\nabla v|^2 + a_0 |v|^2 \, dx - \int_\Omega f_0 v\, dx \right \}, 
\end{equation}
The theory of Lagrange multiplier can allow to reformulate the equation \eqref{eq:orig} as the following minmax problem \cite{babuvska1973finite}, which finds $u \in H^1(\Omega)$ and $\lambda \in H^{-\frac{1}{2}}(\Gamma)$ that minimize the following minmax optimization problem:   
\begin{equation}\label{eq:saddle}  
\min_{v \in H^1(\Omega)} \left \{ \max_{\eta \in H^{-\frac{1}{2}}(\Gamma)} \frac{1}{2}\int_\Omega |\nabla v|^2 + a_0|v|^2 \, dx - \int_\Gamma \eta (v - g) \, ds - \int_\Omega f_0 v\, dx \right \},  
\end{equation}
where $\lambda$ is the Lagrange multiplier that ensures the boundary condition that $u = g$ on $\partial \Omega$.

To discuss the operator form of the equation for \eqref{eq:saddle}, we shall introduce two operators $A_0 : H^{1}(\Omega) \mapsto H^{-1}(\Omega)$ and $B : H^1(\Omega) \mapsto H^{\frac{1}{2}}(\Gamma)$ defined through: 
\begin{subeqnarray}
\langle A_0 u , v\rangle &=& \int_\Omega \nabla u \cdot \nabla v + a_0 uv \, dx , \quad \forall u, v \in H^1(\Omega)  \\ 
\langle B u , q\rangle 
&=& \int_\Gamma {\rm tr}(u) q \, ds = \int_\Gamma uq \, ds, \quad \forall u \in H^1(\Omega) \mbox{ and } \forall q \in H^{-\frac{1}{2}}(\Gamma), 
\end{subeqnarray}
where ${\rm tr} : H^1(\Omega) \mapsto H^{\frac{1}{2}}(\Gamma)$ is the standard trace operator. The minmax problem \eqref{eq:saddle} can then be written as the following saddle point problem: 
\begin{subeqnarray}\label{eqmain} 
\langle A_0 u, v \rangle + \langle \lambda, B v \rangle &=& \langle f_0, v\rangle, \quad \forall v \in H^1(\Omega) \\ 
\langle Bu, \mu\rangle &=& \langle g, \mu \rangle, \quad \forall \mu \in H^{-\frac{1}{2}}(\Gamma)
\end{subeqnarray}
or simply as the following operator form of equation \eqref{eq:saddle}: for $f_0 \in H^{-1}(\Omega)$ and $g \in H^{\frac{1}{2}}(\Gamma)$, find $u \in H^1(\Omega)$ and $\lambda \in H^{-\frac{1}{2}}(\Gamma)$ such that  
\begin{equation}\label{saddle:eq} 
\begin{pmatrix} 
A_0 & B^* \\ 
B & 0 
\end{pmatrix} 
\begin{pmatrix}
u \\ \lambda 
\end{pmatrix}
= \begin{pmatrix}
f_0 \\ g 
\end{pmatrix}
\end{equation}
Note that for $a_0 > 0$, the operator $A_0$ is isormorphic. On the other hand, for $a_0 = 0$, 
\begin{equation}\label{null} 
{\rm ker}(A) = {\rm span}\{1\} \quad \mbox{ and } \quad {\rm ker}(B) = H^1_0(\Omega).  
\end{equation}
Thus the induced norm by the operator $A_0$ is seminorm on $H^1(\Omega)$. 
We remark that $\lambda + \nabla u \cdot n = 0$ on $\Gamma$. Due to the equation \eqref{null}, we note that the problem \eqref{saddle:eq} can be shown to be true by establishing the following inf-sup condition
\begin{equation}\label{infsupc} 
\inf_{q \in H^{-\frac{1}{2}}(\Gamma)} \sup_{u \in H^1(\Omega)} \frac{\langle q, Bu \rangle}{\|u\|_1} \gtrsim \|q\|_{-\frac{1}{2}}. 
\end{equation}
This inf-sup condition is actually a consequence of the following Theorem:  
\begin{Theorem}[Theorem 2.7 \cite{babuvska1973finite}]\label{infsupp} 
Let $q \in H^{-\frac{1}{2}}(\Gamma)$ and $u$ be the weak solution of the following problem: 
\begin{subeqnarray}\label{weakp}  
-\Delta u + u &=& 0 \quad \mbox{ in } \Omega \\ 
\nabla u \cdot n &=& q \quad \mbox{ on } \Gamma, 
\end{subeqnarray}
then the following inequalities hold: 
\begin{equation}\label{first} 
\|q\|^2_{-\frac{1}{2}} \lesssim \int_\Gamma q u \, ds = \|u\|_1^2.  
\end{equation}
\end{Theorem}
\begin{proof}
Since $u$ solves the equation \eqref{weakp}, we have that 
\begin{equation}
\|u\|_1^2 = \int_\Omega |\nabla u|^2 + |u|^2 \, dx = \int_\Gamma q u \, ds. 
\end{equation}
Thus, we only need to show that the first inequality in the equation \eqref{first}. For any given $v \in H^{\frac{1}{2}}(\Gamma)$, let $w$ be the weak solution to the following equation:
\begin{subeqnarray*}
-\Delta w + w &=& 0, \quad \mbox{ in } \Omega \\ 
w &=& v, \quad \mbox{ on } \Gamma. 
\end{subeqnarray*}
Then, it holds that $\|w\|_1 \lesssim \|v\|_{\frac{1}{2}}$. On the other hand, we consider the following inequality: 
\begin{equation}
\frac{\left | \int_\Gamma q v \, ds \right |}{\|v\|_{\frac{1}{2}}}  = \frac{\left | \int_{\Omega} \nabla u \cdot \nabla w + u w \, dx \right |}{\|v\|_{\frac{1}{2}}} \lesssim \frac{\|u\|_{1} \|w\|_1}{\|v\|_{\frac{1}{2}}} \lesssim \|u\|_1.  
\end{equation}
Since $v$ is arbitrary, we arrive at 
\begin{equation}
\|q\|_{-\frac{1}{2}} \lesssim \|u\|_1. 
\end{equation}
This completes the proof. 
\end{proof}
This leads to the following inf-sup condition. 
\begin{Theorem}\label{infsupx}
The following inf-sup condition holds: 
\begin{equation}\label{infsup} 
\inf_{q \in H^{-\frac{1}{2}}(\Gamma)} \sup_{u \in H^1(\Omega)} \frac{\langle q, Bu \rangle}{\|u\|_1} \gtrsim \|q\|_{-\frac{1}{2}}. 
\end{equation}
\end{Theorem} 
\begin{proof} 
Due to Theorem \ref{infsupp}, we have that for any given $q \in H^{-\frac{1}{2}}(\Gamma)$, there exists $u \in H^1(\Omega)$ such that 
\begin{equation}
\|u\|_{1} \|q\|_{-\frac{1}{2}} \lesssim \langle q, Bu \rangle. 
\end{equation}
Thus, we have that 
\begin{equation}
\sup_{u \in H^1(\Omega)} \frac{\langle q, B u \rangle}{\|u\|_1} \gtrsim \|q\|_{-\frac{1}{2}}. 
\end{equation}
This completes the proof. 
\end{proof}
In passing to the next subsection, we shall identify the space $H^{-\frac{1}{2}}(\Gamma)$ and $H^{\frac{1}{2}}(\Gamma)$, using the Riesz representation theorem as usual \cite{goodrich1970riesz,xu2003some}. We also shall then consider the equivalent formulation to \eqref{saddle:eq}: for $f \in L^2(\Omega)$ and $g \in H^{\frac{1}{2}}(\Gamma)$, find $u \in H^1(\Omega)$ and $\lambda \in H^{\frac{1}{2}}(\Gamma)$ such that  
\begin{equation}\label{saddle:eqx} 
\begin{pmatrix} 
A & B^* \\ 
B & 0 
\end{pmatrix} 
\begin{pmatrix}
u \\ \lambda 
\end{pmatrix}
= \begin{pmatrix}
f \\ g 
\end{pmatrix}, 
\end{equation}
where 
\begin{equation}
\langle Au,v\rangle = \langle A_0 u, v\rangle + (u,v)_{0,\Gamma} \quad \mbox{ and } \quad f = f_0 + B^* g,  
\end{equation}
where $(\cdot,\cdot)_{0,\Gamma}$ is the $L^2$ inner product on $\Gamma$. Note that due to the Riesz representation theorem, we can see $B^*$ as a mapping from $H^{\frac{1}{2}}(\Gamma)$ to $H^{-1}(\Omega)$ (see also \cite{bacuta2006unified}). We note that the main reason to consider the equation \eqref{saddle:eqx} is to ensure the invertibility of $A$. Note that the operator $A$ is isomorhic from $H^1(\Omega)$ to $H^{-1}(\Omega)$ and induces the norm, $\langle Au,u\rangle = \|u\|_A^2$. Furthermore, the following Lemma holds \cite{evans2022partial}: 
\begin{lemma}
It holds true that 
\begin{equation} 
\|u\|_1^2 \lesssim \|u\|_A^2 \lesssim \|u\|_1^2, \quad \forall u \in H^1(\Omega).  
\end{equation}
\end{lemma} 
Thus, we shall not distinguish two norms $\|\cdot\|_1$ and $\|\cdot\|_A$. 

\subsection{The Augmented Lagrangian Uzawa Method in the PDE Level} 

Typically, the augmented Lagrangian method has been applied to the discrete system. Namely, once two discrete spaces for $u$ and $\lambda$ are given, say $U_n$ and $W_n$, both of which are a class of shallow neural network functions, the discrete system can then be formulated as follows: 
\begin{eqnarray}\label{mixed} 
\langle A u_n, v_n \rangle + \langle \lambda_n, B v_n \rangle = \langle f,v_n \rangle, \quad \forall v_n \in U_n \subset H^1(\Omega) \\ 
\langle B u_n, \mu_n \rangle = \langle g, \mu_n \rangle, \quad \forall \mu_n \in W_n \subset H^{\frac{1}{2}}(\Gamma).  
\end{eqnarray}
The main issue in the discrete formulation \eqref{mixed} lies at the need of the discrete version of the inf-sup condition, i.e., 
\begin{equation}
\sup_{v_n \in U_n} \frac{\langle \lambda_n, B v_n\rangle }{\|v_n\|_1} \gtrsim \|\mu\|_{0,\Gamma}, \quad \forall \mu_n \in W_n.     
\end{equation}
Generally, the establishment of the discrete inf-sup condition is nontrivial. On the other hand, the discrete inf-sup condition leads to the convergence of the solution $(u_n,\lambda_n)$ to $(u,\lambda)$. Furthermore, the construction of $u_n$ and $\lambda_n$ can be made via the well-known augmented Lagrangian Uzawa method in a speedy way \cite{xu1992iterative,leethesis2004,lee2007robust}. However, the discrete inf-sup condition is in general immature in the setting of the shallow neural network approximations to the best knowledge of the authors. In fact, this is expected to stay as the major challenge for the shallow neural network spaces as they are nonlinearly constructed spaces. As such the weak imposition of the Dirichlet boundary condition is handled via introducing the penalty term in the PDEs \cite{cai2020deep,xu2020finite,yu2018deep}. The main advantage in this approach is that the penalized PDEs are elliptic and thus, the choice of approximate spaces for the approximate solution does not require the limitation of the inf-sup condition. Furthermore, the resulting discrete system can be rather easier to deal with when training algorithms are applied. Such techniques have in fact been successful and popular. The notable difference here lies in that the aformentioned penalty methods are formulated in the continuous level while the penalty methods have been used popularly to handle the weakly imposed symmetry in formulations of the discrete PDE levels, which include Nitsche's technique (see \cite{nitsche1972dirichlet,stenberg1995some} and references cited therein).  

In this section, we shall introduce the augmented Lagrangian Uzawa method constructed in the purely continuous level with the specified penalty parameter, denoted by $\delta \ll 1$. This leads us to obtain a sequence of penalized PDEs and corresponding solutions $\{(u_\delta^\ell ,\lambda_\delta^\ell)\}_{\ell=1,2,\cdots}$, which converge to the original PDE with the sequence of solutions convergent to the solution $(u,\lambda)$ to the original PDE. Namely, we shall show that for $\delta > 0$, it holds that 
\begin{equation}
\lim_{\ell \rightarrow \infty} u_\delta^\ell = u \quad \mbox{ and } \quad \lim_{\ell \rightarrow \infty} \lambda_\delta^\ell = \lambda.
\end{equation}
Our observation is that the first penalized PDE in the sequence of PDEs is the penalized PDE discussed in 
\cite{cai2020deep,xu2020finite,yu2018deep} (see also \cite{girault2012finite}) under the condition that the initial Lagrange multiplier is chosen as zero. The issue arises here that if such a first penalized PDE is solved numerically to obtain the numerical solution $u_n$ using the shallow neural network approximations, then the error estimate in $H^1$ norm may not be sufficiently good. In fact, our numerical experiments demonstrate that there is a difficulty to obtain an optimal error estimate. The main contribution in this paper lies at suggesting to use the second iterate among the sequence of penalized PDEs. Theoretically, we have been able to provide an optimal error estimate for the numerical solution to the second step penalized PDE. 

For a given $\delta > 0$, we consider the following stabilized version of the PDE \eqref{eqmain}, which is given as follows: 
\begin{equation}\label{eqmainx}
\begin{pmatrix} 
A + \delta^{-1} B^*B & B^* \\ 
B & 0 
\end{pmatrix} \begin{pmatrix} u \\ \lambda \end{pmatrix} = \begin{pmatrix} f + \delta^{-1} B^*g \\ g \end{pmatrix}. 
\end{equation}
It is evident that this problem \eqref{eqmainx} is equivalent to \eqref{eqmain}. We now consider to generate a sequence of penalized PDEs which approximate the original PDE. This can be done by applying the Uzawa iterative method. Namely, we apply the augmented Lagrangian Uzawa iteration to obtain the following iterates: Given $(u_\delta^\ell, \lambda_\delta^\ell)$, the new iterate $(u_\delta^{\ell+1}, \lambda_\delta^{\ell+1})$ is obtained by solving the following system:
\begin{subeqnarray}\label{uzawa} 
(A + \delta^{-1} B^*B) u_\delta^{\ell+1} &=& f - B^* \lambda_\delta^\ell + \delta^{-1} B^*g, \\ 
\lambda_\delta^{\ell+1} &=& \lambda_\delta^\ell + \delta^{-1} (B u_\delta^{\ell+1} - g). 
\end{subeqnarray}
One can view the aforementioned augmented Lagrangian Uzawa iterative method as the generation of penalized PDEs that can lead to the original equation. Note that the sequence of penalized PDEs are decoupled elliptic equations for $u_\delta^\ell$ and $\lambda_\delta^\ell$. Using the inf-sup condition in Therem \ref{infsupx}, we can establish the following error bounds in the purely continuous level:
\begin{Theorem}\label{main:thm1} 
The iterate $(u_\delta^\ell,\lambda_\delta^\ell)$ for $\ell = 0,1,\cdots$, satisfies the following estimate: for $\ell \geq 0$,  
\begin{subeqnarray*}
\|\lambda - \lambda_\delta^\ell\|_{0,\Gamma} &\leq& \left ( \frac{\delta}{\delta + \alpha_0} \right )^\ell \|\lambda - \lambda_\delta^0\|_{0,\Gamma}, \\ 
\|u - u_\delta^{\ell+1}\|_A &\leq& \sqrt{\delta} \|\lambda - \lambda_\delta^{\ell}\|_{0,\Gamma} \leq \sqrt{\delta} \left ( \frac{\delta}{\delta + \alpha_0} \right )^\ell \|\lambda - \lambda_\delta^0\|_{0,\Gamma},  
\end{subeqnarray*} 
where 
\begin{equation}
\alpha_0 = \inf_{\mu \in H^{\frac{1}{2}}(\Gamma)}\sup_{v \in H^1(\Omega)} \frac{\langle \mu, Bv \rangle^2}{\langle Av,v\rangle^2}.  
\end{equation} 
\end{Theorem}
\begin{proof} 
It is clear that $A + \delta^{-1} B^*B : H^1(\Omega) \mapsto H^{-1}(\Omega)$ is isomorphic. We define the Schur complement operator $S_\delta$ defined by 
\begin{equation}
S_\delta = B(A + \delta^{-1} B^*B)^{-1} B^*. 
\end{equation}
The application of the Sherman-Morrison-Woodbury formula (see \cite{ojeda2020fast,deng2011generalization} and references cited therein) leads to 
\begin{equation}
S_\delta^{-1} = (B(A + \delta^{-1} B^* B)^{-1} B^*)^{-1} = \delta^{-1}I + (B A^{-1} B^*)^{-1}. 
\end{equation}
We then note that the straighforward calculation (see \cite{bramble1997iterative} and \cite{bacuta2006unified}) leads that 
\begin{equation}
(B A^{-1} B^*\mu,\mu) = \sup_{v \in H^1(\Omega)} \frac{\langle \mu, Bv\rangle^2}{\|v\|_A^2} 
\end{equation}
Thus, the operator $S = B A^{-1}B^*$ is well-conditioned due to the inf-sup condition \eqref{infsup}. At this stage, using the argument presented in \cite{ojeda2020fast,lee2007robust,leethesis2004}, we obtain the estimates. This completes the proof. 
\end{proof}
The main advantage in the formulation given as the augmented Lagrangian Uzawa \eqref{uzawa} is that the sequence of equations to be solved are decoupled elliptic PDEs. For $\ell = 1,2,\cdots,$,  $u_{\delta}^{\ell}$ and $\lambda_\delta^\ell$ solve the following decoupled PDEs: for all $v \in H^1(\Omega)$, 
\begin{subeqnarray}
a_\delta(u_{\delta}^\ell,v) &=& (f,v) + \delta^{-1}(g, v)_{0,\Gamma} - (\lambda_\delta^{\ell-1},v)_{0,\Gamma}, \\ 
(\lambda_\delta^{\ell}, \mu)_{0,\Gamma} &=& (\lambda_\delta^{\ell-1},\mu)_{0,\Gamma} + \delta^{-1}(u_{\delta}^{\ell}, \mu)_{0,\Gamma} - (g, \mu)_{0,\Gamma}, \quad \forall \mu \in H^{\frac{1}{2}}(\Gamma), 
\end{subeqnarray}
where 
\begin{equation}
a_\delta(u,v) := a(u,v) + \delta^{-1}(u,v)_\Gamma, \quad \forall u,v \in H^1(\Omega),   
\end{equation}
with $a(u,v) = \langle Au,v\rangle, \quad\forall u, v \in H^1(\Omega)$. We note that the inner produce $a_\delta(\cdot,\cdot) : H^1(\Omega) \times H^1(\Omega) \mapsto \Reals{}$ induces the norm on $H^1(\Omega)$, which is denoted by $\|\cdot\|_{a_\delta}$. Note that the following estimate holds: 
\begin{equation}
\|u - u_{\delta}^\ell\|_A = O \left (\delta^{\ell - \frac{1}{2}}\right ).  
\end{equation}
We observe that if $\lambda_\delta^0$ is chosen to be zero, then the first iterate of $u_\delta^1$ can still make sense as a good approximate solution to $u$ and $u_\delta^1$ is the solution to the following PDE:   
\begin{equation}
(A + \delta^{-1} B^*B ) u_\delta^1 = f + \delta^{-1} B^* g   
\end{equation} 
We remark that in the pioneering works \cite{yu2018deep, xu2020finite, cai2020deep}, the elliptic equations with the Dirichlet boundary condition are proposed to be solved as the following minimization problem: 
\begin{equation}\label{penalty} 
u_\delta = {\rm arg} \min_{v \in H^1(\Omega)} \left \{ \frac{1}{2} a(v,v) + \delta^{-1} \int_\Gamma |v - g|^2 \, ds - \int_\Omega fv\, dx \right \}. 
\end{equation} 
Here $\delta$ is introduced as a penalty paramter. In fact, such an approach is very much popular in general. The formulation leads to a solution $u_\delta$ for an appropriate choice of $\delta$ in numerical experiments. We observe that it is nothing else than the first iterate of the augmented Lagrangian Uzawa method. 
\begin{lemma}
The solution $u_\delta$ to \eqref{penalty} is equivalent to $u_\delta^1$ when $\lambda_\delta^0 = 0$ in the algorithm \eqref{uzawa} and it holds that 
\begin{equation}\label{delta}
\|u - u_\delta\|_1 \leq \sqrt{\delta} \|\lambda\|_{0,\Gamma}.  
\end{equation}
\end{lemma}
\begin{proof} 
The equivalence between $u_\delta^1$ and $u_\delta$ is evident and the error esimate \eqref{delta} is from Theorem  \ref{main:thm1}. 
This completes the proof. 
\end{proof} 
Namely, the finite neuron method has used the case when $\delta = 1$ and the supercloseness has been justified with $\sqrt{\delta}$ in \cite{xu2020finite}. This is not new. A similar discussion was made in \cite{girault2012finite} to solve Stokes equation and also in   \cite{adhikari;kim;lee;sheen2023} to solve Darcy's flow. Again, the main advantage of using the intermediate function $u_\delta$ is that it is the solution to the elliptic PDEs and thus, the standard solution technique like orthogonal greedy algorithm can be used to obtain an optimal solution $u_n$ that approximates $u_\delta$. The solution $u_n$ is taken as an approximate solution to $u$. The main issue in this approach is to involve $u_\delta$ in the error analysis of $u_n$ against $u$. It is expected that the closer the $u_\delta$ to $u$ is, the smaller the difference between $u_n$ and $u$ is. However, the use of $u_\delta$ as an intermediate solution to obtain $u$ introduces a significant difficulty to establish the optimal error \cite{xu2020finite} due to the penalty paramter $\delta^{-1}$. The following estimate was the state-of-the art result established in \cite{xu2020finite}: 
\begin{Theorem}
Let $u_n$ be the shallow neural network approximation of $u_\delta^1$ and $\delta \approx n^{-\frac{1}{2} - \frac{2(k-1)+1}{2d}}$. Then we have that 
\begin{eqnarray*}
\|u - u_n\|_1 \leq n^{-\frac{1}{4} - \frac{2(k-1)+1}{4d}}. 
\end{eqnarray*}
\end{Theorem} 
Thus, the optimal convergence rate could not be achieved. We note that numerically, it was shown that the optimal error bound can be obtained under the condition that $\delta$ is chosen sufficiently small. We remark that the use of $\delta$, too small can lead to illconditioned system to be solved, which is an additional bottleneck in computations.  

A main idea in our paper is not to use $u_\delta^1$ as an intermediate solution to obtain the approximate solution to $u$, but rather to use $u_\delta^2$, i.e., the second iterate of the augmented Lagrangian Uzawa method \eqref{uzawa}. The use of $u_\delta^2$ as an intermediate solution can lead to the following estimate: 
\begin{equation}
\|u - u_n\|_1 \lesssim \delta^{3/2} + n^{-\frac{1}{2} - \frac{2(k-1)+1}{2d}},  
\end{equation}
which will be established in \S \ref{opt}. This means that $\delta$ can be chosen to behave like $\left [ n^{-\frac{1}{2} - \frac{2(k-1)+1}{2d}}\right ]^{2/3}$ to get the optimal error estimate. Namely, the choice of $\delta$ is relatively large. Numerically, we show that such a choice of $\delta$ does not lead to the optimal accuracy when $u_\delta^1$ is chosen as an intermediate solution in general as discussed in \S \ref{hahi}. 

\section{Optimality in $H^1$ norm of the shallow neural network solution for $u$}\label{opt} 
In this section, we shall discuss the shallow neural network solution $u_n$ to $u_{\delta}^{2}$ and its optimal accuracy to approximate $u$. We note that to obtain the approximate solution to $u_\delta^2$, a couple of steps are needed. Among others, one has to solve for $u_\delta^1$ numerically. After this is done, $\lambda_\delta^1$ is approximated. Lastly, we are in a position to obtain the numerical solution to approximate $u_{\delta}^2$. This is summarized in the Algorithm \ref{aug}. 

\begin{algorithm}[h]  
\caption{}\label{aug}
Given $k$ and $\delta = \left [ n^{-1/2 - \frac{2(k-1)+1}{2d}} \right ]^{2/3}$,  
\begin{algorithmic}[1]
\STATE{Solve for $u_{\delta,n}^1 \in \Sigma_d^n(\sigma)$
\begin{equation}
a(u_{\delta,n}^{1},v) + \delta^{-1}(u_{\delta}^{1},v)_{\Gamma} = \langle f, v \rangle + \delta^{-1}(g, v)_{0,\Gamma}, \quad \forall v \in \Sigma_d^n(\sigma),    
\end{equation}
for $u_\delta^1$ in a shallow neural network space, which is denoted by $u_{\delta,n}^1$}
\STATE{Define $\lambda_n \in \Sigma_{d-1}^n(\sigma)$ by the following identity:  
\begin{equation}
\lambda_n := \delta^{-1} \left (B u_{\delta,n}^1 - g \right ).
\end{equation}}
\STATE{Solve for $u_n$ in $\Sigma_d^n(\sigma)$
\begin{equation}
a(u_n,v) + \delta^{-1}(u_n,v)_{0,\Gamma} = \langle f, v \rangle + \delta^{-1} (g, v)_{0,\Gamma} - (\lambda_{n}, v)_{0,\Gamma}, \quad \forall v \in H^1(\Omega). 
\end{equation}
}
\end{algorithmic}
\end{algorithm}
For the sake of the theoretical analysis, we shall introduce an additional numerical solution, denoted by $u_{\delta,n}^2 \in \Sigma_d^n(\sigma)$, which is defined as follows: 
\begin{equation}\label{approx:eq}
a_{\delta}(u_{\delta,n}^2,v) = \langle f, v \rangle + \delta^{-1}(g, v)_{0,\Gamma} - ( \lambda_{\delta}^1,v)_{0,\Gamma}, \quad \forall v \in \Sigma_{d}^n(\sigma). 
\end{equation} 

The following trace theorem is the instrumental to obtain the error estimate. 
\begin{lemma}\label{trace} 
For any $\epsilon > 0$, the following trace theorem holds: 
\begin{equation}
\|\phi\|_{0,\Gamma}^2 \lesssim \epsilon^{-1} \|\phi\|_{0}^2 + \epsilon |\phi|_{1}^2, \quad \forall \phi \in H^1(\Omega). 
\end{equation}
\end{lemma}
We shall first establish the following instrumental lemma. The idea lies in the duality argument and universal approximation principle of the shallow neural network functions. 
\begin{lemma}\label{main:tool} 
The folloing bound holds true: 
\begin{equation}
\|u_\delta^2 - {u}_{\delta,n}^2\|_0 \lesssim n^{-\frac{1}{2} - \frac{2(k-1)+1}{2d}}  (1 + \delta^{-\frac{1}{2}} ) 
 \|u_\delta^2 - {u}_{\delta,n}^2\|_{a_\delta}.  
\end{equation}
\end{lemma}
\begin{proof} 
First, we consider the definition of the $L^2$ norm: 
\begin{equation}
\|u_\delta^2 - {u}_{\delta,n}^2\|_0 = \sup_{\psi \in L^2(\Omega), \|\psi\|_0 = 1} (u_\delta^2 - {u}_{\delta,n}^2, \psi)_0. 
\end{equation}
We shall now invoke the Aubin-Nitsche duality argument \cite{schatz1996some} that for $\psi \in L^2$ with $\|\psi\|_0 = 1$, let $v \in H^1(\Omega)$ be the solution to the following equation: 
\begin{equation}
a_\delta(\eta,v) = (\psi,\eta)_0, \quad \forall \eta \in H^1(\Omega). 
\end{equation}
Then, we have that 
\begin{eqnarray*}
|(u_\delta^2 - {u}_{\delta,n}^2, \psi)| &=& |a_{\delta}(u_\delta^2 - {u}_{\delta,n}^2, v)| \\
&=& |a_\delta(u_\delta^2 - {u}_{\delta,n}^2, v - \chi)| \\  &\lesssim& \|u_\delta^2 - {u}_{\delta,n}^2\|_{a_\delta} \|v - \chi\|_{a_\delta},   
\end{eqnarray*}
where $\chi$ is any function in $\Sigma_d^n(\sigma)$. On the other hand, we have that 
\begin{equation}
\|v - \chi\|_{a_\delta} \lesssim (1 + \delta^{-\frac{1}{2}}) \|v - \chi\|_1.  
\end{equation}
Thus, we arrive at the desired estimate:
\begin{eqnarray*}
\|u_\delta^2 - {u}_{\delta,n}^2\|_0 &\lesssim& (1 + \delta^{-\frac{1}{2}}) \inf_{\chi \in \Sigma_{d}^n(\sigma)} \|v - \chi\|_1 \|u_\delta^2 - {u}_{\delta,n}^2\|_{a_\delta} \\
&\lesssim& n^{-\frac{1}{2} - \frac{2(k-1)+1}{2d}}  (1 + \delta^{-\frac{1}{2}} ) 
 \|u_\delta^2 - {u}_{\delta,n}^2\|_{a_\delta}.  
\end{eqnarray*}
This completes the proof. 
\end{proof} 
We are now in a position to establish the error estimate for ${u}_{\delta,n}^2$ against $u$. 
\begin{Theorem}\label{main:thm} 
Let $\delta$ be chosen as follows:
\begin{equation}
\delta \lesssim n^{-\frac{1}{3} - \frac{2(k-1)+1}{3d}}. 
\end{equation}
Then the shallow neural network solution ${u}_{\delta,n}^2$ for the equation \eqref{approx:eq} satisfies the following error bound:   
\begin{equation}
\|u_\delta^2 - {u}_{\delta,n}^2\|_1 \lesssim n^{-\frac{1}{2} - \frac{2(k-1)+1}{2d}}.   
\end{equation}
\end{Theorem}
\begin{proof}
Since ${u}_{\delta,n}$ solves \eqref{approx:eq}, we have that 
\begin{equation}
\|u_\delta^2 - {u}_{\delta,n}^2\|_{a_\delta} \leq \|u_\delta^2 - {u}_{\delta,n}^2\|_{1} + \delta^{-\frac{1}{2}} \|u_\delta^2 - {u}_{\delta,n}^2\|_{0,\Gamma}
\end{equation}
Due to Lemma \ref{trace}, we have that for $c$ being a sufficiently large generic constant, 
\begin{equation} 
\delta^{-\frac{1}{2}} \|u_\delta^2 - {u}_{\delta,n}^2\|_{0,\Gamma} \lesssim \delta^{-\frac{1}{2}} \left ( (c\delta)^{-\frac{1}{2}} \|u_\delta^2 - {u}_{\delta,n}^2\|_{0} + (c\delta)^{\frac{1}{2}} |u_\delta^2 - {u}_{\delta,n}^2|_{1} \right ). 
\end{equation} 
Thus, we arrive at 
\begin{equation} 
\|u_\delta^2 - {u}_{\delta,n}^2\|_{a_\delta} \lesssim \|u_\delta^2 - {u}_{\delta,n}^2 \|_{1} + c^{-\frac{1}{2}} \delta^{-1} \|u_\delta^2 - {u}_{\delta,n}^2\|_{0} + c^{\frac{1}{2}} |u_\delta^2 - {u}_{\delta,n}^2|_{1}
\end{equation} 
We now apply the Lemma  \ref{main:tool}, i.e.,
\begin{equation}
\|u_\delta^2 - {u}_{\delta,n}^2\|_0 \lesssim \left ( \delta^{\frac{3}{2}} + \delta \right ) \|u_\delta^2 - {u}_{\delta,n}^2\|_{a_\delta}
\end{equation}
to get that 
\begin{subeqnarray} 
\|u_\delta^2 - {u}_{\delta,n}^2\|_{a_\delta} &\lesssim& \|u_\delta^2 - {u}_{\delta,n}^2\|_{1} + c^{-\frac{1}{2}} \|u_\delta^2 - {u}_{\delta,n}^2\|_{a_\delta}.  
\end{subeqnarray}
The Kick-Back argument leads to the following bound: 
\begin{equation}
\|u_\delta^2 - {u}_{\delta,n}^2\|_{a_\delta} \lesssim \|u_\delta^2 - {u}_{\delta,n}^2\|_1 \lesssim n^{-\frac{1}{2} - \frac{2(k-1)+1}{2d}}. 
\end{equation}
This completes the proof. 
\end{proof} 
We note that Algorithm \ref{aug} requires to obtain the solution $u_n$ that approximate the equation given as follows: 
\begin{equation}\label{approx:eq1}
a_\delta(\overline{u}_\delta^2,v) = \langle f,v \rangle + \delta^{-1}(g, v)_{0,\Gamma} - (\lambda_{n}, v)_{0,\Gamma}, \quad \forall v \in H^1(\Omega). 
\end{equation} 
Basically the argument used in Lemma  \ref{main:tool} and Theorem \ref{main:thm}, leads to obtain the following corollary: 
\begin{corollary}
Let $\delta$ be chosen as follows:
\begin{equation}
\delta \lesssim n^{-\frac{1}{3} - \frac{2(k-1)+1}{3d}}. 
\end{equation}
Then the shallow neural network solution ${u}_{n}$ for the equation \eqref{approx:eq1} satisfies the following error bound:   
\begin{equation}
\|\overline{u}_\delta^2 - {u}_n\|_1 \lesssim n^{-\frac{1}{2} - \frac{2(k-1)+1}{2d}}.   
\end{equation}
\end{corollary}
We are now in a position to prove the main result in this paper. 
\begin{Theorem}\label{main:thm2} 
If $\lambda_n$ is obtained so that the following holds true:
\begin{equation}
\|\lambda_\delta^1 - \lambda_n\|_{\frac{1}{2},\Gamma} \lesssim n^{-\frac{1}{2} - \frac{2(k-1) + 1}{2d}}, 
\end{equation}
Then we get 
\begin{equation}
\|u - u_n\|_1 \lesssim n^{-\frac{1}{2} - \frac{2(k-1) + 1}{2d}}. 
\end{equation}
\end{Theorem}
\begin{proof} 
A simple triangle inequality leads to 
\begin{equation}
\|u - u_n\|_1 \leq \|u - u_\delta^2\|_1 + \|u_\delta^2 - \overline{u}_{\delta}^2\|_1 + \|\overline{u}_\delta^2 - u_n\|_1. 
\end{equation}
The first and third terms are bounded by $n^{-\frac{1}{2} - \frac{2(k-1)+1}{2d}}$ under the specified choice of $\delta$. Thus, it remains to estimate the term in the middle. 
\begin{eqnarray*}
\|u_{\delta}^{2} - \overline{u}_\delta^2\|_1^2 &\leq& \|u_\delta^2 - \overline{u}_\delta^2\|_{a_\delta}^2 \\
&\leq& (\lambda_{\delta}^1 - \lambda_n, u_{\delta}^{2} - \overline{u}_{\delta}^{2})_{0,\Gamma} \\
&\lesssim& \|\lambda_\delta^1 - \lambda_n\|_{0,\Gamma} \|u_\delta^2 - \overline{u}_{\delta}^{2}\|_1. 
\end{eqnarray*}
Thus, we have the following estimate:
\begin{eqnarray} 
\|u_\delta^2 - \overline{u}_\delta^2\|_1 \lesssim 
\|\lambda_\delta^1 - \lambda_n\|_{0,\Gamma}. 
\end{eqnarray}
This completes the proof. 
\end{proof} 
Theorem \ref{main:thm2} indicates that it is desirable to obtain the approximate Lagrange multiplier accurately to obtain the optimal error estimate. This can be in fact achieved if we solve $u_\delta^1$ rather accurately, possibly by choosing larger $k$. On the other hand, the optimal order was observed to hold even if $\lambda_n$ does not satisfy the desired accuracy as presented in \S \ref{num:ex}.  

\section{Orthogonal Greedy Algorithm}
In this section, we now discuss the orthogonal greedy algorithm to obtain the shallow neural network solution to the elliptic PDEs. Given a target function $u$ in a Hilbert space with an inner product $\langle \cdot, \cdot \rangle_H$ and the induced norm $\|\cdot\|_H$, we consider to approximate $u$ in $\|\cdot\|_H$ norm, using a sparse linear combinations of dictionary elements, 
\begin{equation}
u_n = \sum_{i=1}^n a_i g_i,
\end{equation}
where $g_i \in \mathbb{D}$ depends on the function $u$. For the general dictionary $\mathbb{D} \subset H^1(\Omega)$, the orthogonal greedy algorithm, also known as orthogonal matching pursuit is given as the Algorithm \ref{OGA} can be shown to obtain $u_n$ that satisfies the estimate given in Lemma below \cite{siegel2021greedy}.  
\begin{algorithm}[h]
\caption{Orthogonal Greedy Algorithm}\label{OGA}
Given $u_0 = 0$, we perform 
\begin{algorithmic}[1]
\STATE{STEP 1: Solve the following optimization problem:   
\begin{equation}\label{alg1min}
g_k = {\rm arg}\max_{g \in \mathbb{D}} \left |\langle u - u_{k-1}, g \rangle_H \right | 
\end{equation}}
\STATE{STEP 2: Update the function approximation: 
\begin{equation}
\langle u_k,v \rangle_H = \langle u, v\rangle_H, \quad \forall v \in V_k = {\rm span} \{g_1,\cdots,g_k\}. 
\end{equation}
}
\end{algorithmic}
\end{algorithm}
\begin{lemma}
Given a function $u \in H^m(\Omega)$, the orthogonal greedy algorithm can attain  an $n-$term expansion which satisfies
\begin{equation}\label{sharpx}
\|u - u_n\|_{H^m} \lesssim n^{-\frac{1}{2}}, 
\end{equation} 
where 
\begin{equation}
\Sigma_n(\mathbb{D}) = \left \{ \sum_{i=1}^n a_i g_i, a_i \in \Reals{}, g_i \in \mathbb{D} \right \}. 
\end{equation}
is the $n-$term nonlinear dictionary expansions. 
\end{lemma}
Although the approximation rate given in \eqref{sharpx} is sharp for the general dictionaries, the rate can be improved using a stratified sampling argument for compact dictionaries. 
For shallow neural networks with ${\rm ReLU}^k$ activation function, the relavant dictionary is $\Sigma_d^n(\sigma)$. With such a choice of dictionary, it was proven that the orthogonal greedy algorithm attains the $n$ term expansion, $u_n$ that can lead to the following estimate: 
\begin{equation}
\|u - u_n\|_{H^m} \lesssim n^{-\frac{1}{2} - \frac{2(k-m)+1}{2d}}, 
\end{equation}
which is optimal \cite{siegel2021greedy}. 

\section{Numerical Experiments}
\label{num:ex}
In this section, we shall present numerical results that confirm the theory established in the previous sections. For a fixed $\delta$, we define the inner product for the space of shallow neural networks as folows:  
\begin{equation}
\langle u, v\rangle_H := a_\delta(u,v) = a(u,v) + \delta^{-1} (u,v)_{0,\Gamma}. 
\end{equation}
We shall present some numerical tests that confirm the theory. Our aim is to demonstrate the relationship between $n$ and the $L^2$ and $H^1$ norm convergence of the error. 
\subsection{1D results}\label{hahi} 

In this section, we consider 1D poisson equation with dirichlet boundary condition
\begin{subeqnarray}
-u'' + a_0 u &=& f_0 \quad \mbox{ in } \Omega, \\ 
u&=& 0 \quad \mbox{ on } \Gamma. 
\end{subeqnarray}
We have tested cases with $a_0 = 0$ and $a_0 = 1$ with $\Omega = (-1,1)$ and 
\begin{subeqnarray}
f_0(x) &=& \left(\frac{\pi^2}{4}\right)\cos \left(\frac{\pi}{2} x\right), \mbox{ when } a_0 = 0 \\ 
f_0(x) &=& \left(1+\frac{\pi^2}{4}\right)\cos \left(\frac{\pi}{2} x\right), \mbox{ when } a_0 = 1. 
\end{subeqnarray}
For the numerical quadrature, we use a two-point Gaussian quadrature rule to obtain discrete energy functional after the domain is subdivided into a total of 4000 subintervals. For the OGA, we followed the technique described in \cite{siegel2021greedy}. Numerical test results are listed in the Tables \ref{tab_2} and \ref{tab_3}. 

\subsubsection{The case $k = 1$}
In this section, we present the result that is obtained with $k = 1$.

\begin{table}[H]
\centering
\caption{Numerical results of OGA for 1D elliptic equation with Dirichlet boundary condition and $\delta = 256^{-\frac{2}{3}}$. We show two different results for $a_0=0$ and $a_0=1$. The predicted $H^1$ norm error is about $2\delta^{\frac{3}{2}} = $ 7.8e-03.} 
\vspace{0.5cm}
\label{tab_2}
\renewcommand\arraystretch{1.5}
\scalebox{0.72}{
\begin{tabular}{|c|c|c|c|c|c|c|c|c|}
\hline
\multirow{3}{*}{$n$} &
\multicolumn{4}{c|}{$a_0$=0} &
\multicolumn{4}{c|}{$a_0$=1} \\
\cline{2-9} 
&\multicolumn{2}{c|}{The result of $u_\delta^{1}$} &
\multicolumn{2}{c|}{The result of $u_\delta^{2}$} &
\multicolumn{2}{c|}{The result of $u_\delta^{1}$} &
\multicolumn{2}{c|}{The result of $u_\delta^{2}$} \\
\cline{2-9}  
 & $\|u-u_n\|_{H^1}$ & $n^{-1}$ & $\|u-u_n\|_{H^1}$ & $n^{-1}$& $\|u-u_n\|_{H^1}$ & $n^{-1}$ & $\|u-u_n\|_{H^1}$ & $n^{-1}$  \\ \hline
16 &1.000704e-01 & 1.01  & 8.582538e-02 & 1.21 & 9.656833e-02 & 1.05 & 8.547528e-02 & 1.18 \\ \hline
32  & 6.807205e-02 & 0.56   & 4.099041e-02 & 1.07 & 6.211432e-02 & 0.64  & 4.122272e-02 & 1.05 \\ \hline 
64  & 5.853793e-02 & 0.22  & 2.024263e-02 & 1.02 & 5.124574e-02 & 0.28  & 2.042999e-02 & 1.01 \\ \hline 
128 & 5.596750e-02 & 0.06  & 1.005804e-02 & 1.01 & 4.822536e-02 & 0.09  & 1.017396e-02 & 1.01 \\ \hline 
256  & 5.531920e-02 & 0.02 &  5.047939e-03 & 0.99 & 4.745296e-02 & 0.02  & 5.125002e-03 & 0.99 \\ \hline  
\end{tabular}}
\end{table} 

\subsubsection{The case $k = 2$}
In this section, we present the result that is obtained with $k = 2$.

\begin{table}[H]
\centering
\caption{Numerical results of OGA for 1D elliptic equation with Dirichlet boundary condition and $\delta =256^{-\frac{4}{3}}$. We show two different results for $a_0=0$ and $a_0=1$. The predicted $H^1$ norm error is about $2\delta^{\frac{3}{2}} = $3.0518e-05.} 
\vspace{0.5cm}
\label{tab_3}
\renewcommand\arraystretch{1.5}
\scalebox{0.72}{
\begin{tabular}{|c|c|c|c|c|c|c|c|c|}
\hline
\multirow{3}{*}{$n$} &
\multicolumn{4}{c|}{$a_0$=0} &
\multicolumn{4}{c|}{$a_0$=1} \\
\cline{2-9} 
&\multicolumn{2}{c|}{The result of $u_\delta^{1}$} &
\multicolumn{2}{c|}{The result of $u_\delta^{2}$} &
\multicolumn{2}{c|}{The result of $u_\delta^{1}$} &
\multicolumn{2}{c|}{The result of $u_\delta^{2}$} \\
\cline{2-9}  
 & $\|u-u_n\|_{H^1}$ & $n^{-2}$ & $\|u-u_n\|_{H^1}$ & $n^{-2}$& $\|u-u_n\|_{H^1}$ & $n^{-2}$ & $\|u-u_n\|_{H^1}$ & $n^{-2}$  \\ \hline
16 & 3.577753e-03 & 2.06 & 3.245438e-03 & 2.21 & 3.602155e-03 & 2.09  & 3.234727e-03 & 2.16 \\ \hline 
32  & 1.503081e-03 & 1.25  & 6.744694e-04 & 2.27 & 1.341009e-03 & 1.43  & 6.025989e-04 & 2.42 \\ \hline 
64 & 1.376389e-03 & 0.13  & 1.644393e-04 & 2.04 & 1.202158e-03 & 0.16  & 1.525457e-04 & 1.98 \\ \hline 
 
128  & 1.367156e-03 & 0.01 & 4.035486e-05 & 2.03 & 1.192672e-03 & 0.01 & 3.687647e-05 & 2.05 \\ \hline 
256 & 1.366655e-03 & 0.00   & 9.789934e-06 & 2.04 & 1.192120e-03 & 0.00  & 9.405487e-06 & 1.97 \\ \hline 
\end{tabular}}
\end{table} 

\subsection{2D results}
In this section, we consider 2D the model equation given as follows: 
\begin{equation}
-\Delta u + a_0 u = f_0 \quad \mbox{ in } \Omega,
\end{equation}
subject to the boundary condition that $u = 0$ on $\partial \Omega$. Here $\Omega = (0,1)^2$. 
\begin{subeqnarray}
f_0(x) &=& 2 \pi^{2} \sin \left(\pi x \right)\sin \left(\pi y \right), \mbox{ when } a_0 = 0, \\ 
f_0(x) &=& \left(1+2 \pi^{2}\right)\sin \left(\pi x \right)\sin \left(\pi y \right), \mbox{ when } a_0 = 1. 
\end{subeqnarray}
We use a $4$ point Gaussian quadrature rule in each cell, to obtain the discrete energy functional. More detailed computational details can be found at \cite{siegel2021greedy}. For the boundary integral, we select $1000$ points on each subinterval of $\Gamma$ and use two points Gaussian integration. Numerical results are listed in the Tables \ref{tab_4} and \ref{tab_5}.

\subsubsection{The case $k = 1$}
In this section, we test the case when $k = 1$ is used for the activation ${\rm ReLu}^k$ function. Note that the optimal rate is given as 
\begin{equation}  
\|u - u_n\|_1 \approx n^{-\frac{3}{4}}.
\end{equation}
Thus, if we choose $\delta = n^{-\frac{1}{2}}$, then we expect that the optimal convergence can be achieved.

\begin{table}[H]
\centering
\caption{Numerical results of OGA for 2D elliptic equation with Dirichlet boundary condition and $\delta =256^{-\frac{1}{2}}$. We show two different results for $a_0=0$ and $a_0=1$. The predicted $H^1$ norm error is about $2\delta^{\frac{3}{2}} = $3.12e-02.} 
\vspace{0.5cm}
\label{tab_4}
\renewcommand\arraystretch{1.5}
\scalebox{0.72}{
\begin{tabular}{|c|c|c|c|c|c|c|c|c|}
\hline
\multirow{3}{*}{$n$} &
\multicolumn{4}{c|}{$a_0$=0} &
\multicolumn{4}{c|}{$a_0$=1} \\
\cline{2-9} 
&\multicolumn{2}{c|}{The result of $u_\delta^{1}$} &
\multicolumn{2}{c|}{The result of $u_\delta^{2}$} &
\multicolumn{2}{c|}{The result of $u_\delta^{1}$} &
\multicolumn{2}{c|}{The result of $u_\delta^{2}$} \\
\cline{2-9}  
 & $\|u-u_n\|_{H^1}$ & $n^{-\frac{3}{4}}$ & $\|u-u_n\|_{H^1}$ & $n^{-\frac{3}{4}}$& $\|u-u_n\|_{H^1}$ & $n^{-\frac{3}{4}}$ & $\|u-u_n\|_{H^1}$ & $n^{-\frac{3}{4}}$  \\ \hline
  16  & 2.08e+00 & 0.12  & 1.55e+00 & 0.52 & 2.07e+00 & 0.13 &  1.60e+00 & 0.33 \\ \hline
32 & 9.75e-01 & 1.09 & 8.89e-01 & 0.80 & 9.65e-01 & 1.10 &8.25e-01 & 0.95 \\ \hline   
64  & 3.43e-01 & 1.51  & 3.61e-01 & 1.30 & 3.42e-01 & 1.49&3.53e-01 & 1.22 \\ \hline
128  & 1.55e-01 & 1.15   & 1.61e-01 & 1.16 & 1.55e-01 & 1.14 &  1.63e-01 & 1.12 \\ \hline 
256& 8.22e-02 & 0.92  & 8.08e-02 & 0.99 & 8.22e-02 & 0.92 &  8.03e-02 & 1.02 \\ \hline 
\end{tabular}}
\end{table}

\subsubsection{The case $k = 2$}
In this section, we test the case when $k = 2$ is used for the activation ${\rm ReLu}^k$ function. Note that the optimal rate is given as 
\begin{equation}
\|u - u_n\|_1 \approx n^{-\frac{5}{4}}. 
\end{equation}
Thus, if we choose $\delta = n^{-\frac{5}{6}}$, then we expect that the optimal convergence rate can be achieved.

\begin{table}[H]
\centering
\caption{Numerical results of OGA for 2D elliptic equation with Dirichlet boundary condition and $\delta =256^{-\frac{5}{6}}$. We show two different results for $a_0=0$ and $a_0=1$. The predicted $H^1$ norm error is about $2\delta^{\frac{3}{2}} = $1.95e-03.} 
\vspace{0.5cm}
\label{tab_5}
\renewcommand\arraystretch{1.5}
\scalebox{0.72}{
\begin{tabular}{|c|c|c|c|c|c|c|c|c|}
\hline
\multirow{3}{*}{$n$} &
\multicolumn{4}{c|}{$a_0$=0} &
\multicolumn{4}{c|}{$a_0$=1} \\
\cline{2-9} 
&\multicolumn{2}{c|}{The result of $u_\delta^{1}$} &
\multicolumn{2}{c|}{The result of $u_\delta^{2}$} &
\multicolumn{2}{c|}{The result of $u_\delta^{1}$} &
\multicolumn{2}{c|}{The result of $u_\delta^{2}$} \\
\cline{2-9}  
 & $\|u-u_n\|_{H^1}$ & $n^{-\frac{5}{4}}$ & $\|u-u_n\|_{H^1}$ & $n^{-\frac{5}{4}}$& $\|u-u_n\|_{H^1}$ & $n^{-\frac{5}{4}}$ & $\|u-u_n\|_{H^1}$ & $n^{-\frac{5}{4}}$  \\ \hline
  16  & 2.05e+00 & 0.15&  2.03e+00 & 0.16 & 2.04e+00 & 0.16 &  2.05e+00 & 0.15 \\ \hline 
32 & 6.02e-01 & 1.77 & 6.64e-01 & 1.61 & 6.00e-01 & 1.77 &  4.70e-01 & 2.13 \\ \hline
64  & 1.03e-01 & 2.54 &  1.22e-01 & 2.45 & 1.03e-01 & 2.54 &  9.00e-02 & 2.38 \\ \hline  
128  & 2.37e-02 & 2.12 & 2.16e-02 & 2.49 & 2.37e-02 & 2.12 &  2.34e-02 & 1.94 \\ \hline 
256& 6.27e-03 & 1.92 &  6.23e-03 & 1.79 & 6.24e-03 & 1.92 &  6.45e-03 & 1.86 \\ \hline 
\end{tabular}}
\end{table} 

\subsection{3D results} 
In this section, we consider 3D model equations given as follows: 
\begin{equation}
-\Delta u + a_0 u = f_0 \quad \mbox{ in } \Omega,
\end{equation}
subject to the boundary condition that $u = 0$ on $\partial \Omega$. Here $\Omega = (0,1)^3$. 
\begin{subeqnarray}
f_0(x) &=& 3 \pi^{2} \sin \left(\pi x \right)\sin \left(\pi y \right)\sin \left(\pi z \right), \mbox{ when } a_0 = 0, \\ 
f_0(x) &=& \left(1+3 \pi^{2}\right)\sin \left(\pi x \right)\sin \left(\pi y \right)\sin \left(\pi z \right), \mbox{ when } a_0 = 1. 
\end{subeqnarray}
We use a $8$ point Gaussian quadrature rule in each cell, to obtain the discrete energy functional. More detailed computational details can be found at \cite{siegel2021greedy}. For the boundary integral, we select $1000$ points on each subface of $\Gamma$ and use 4 points Gaussian integration. Numerical results are listed in the Tables \ref{tab_6} and \ref{tab_7}.

\subsubsection{The case $k = 1$}
In this section, we test the case when $k = 1$ is used for the activation ${\rm ReLu}^k$ function. Note that the optimal rate is given as 
\begin{equation}
\|u - u_n\|_1 \approx n^{-\frac{2}{3}}.
\end{equation}
Thus, if we choose $\delta = n^{-\frac{4}{9}}$, then we expect that the optimal convergence rate can be achieved.

\begin{table}[H]
\centering
\caption{Numerical results of OGA for 3D elliptic equation with Dirichlet boundary condition and $\delta =256^{-\frac{4}{9}}$. We show two different results for $a_0=0$ and $a_0=1$. The predicted $H^1$ norm error is about $2\delta^{\frac{3}{2}} = $4.96e-02.} 
\vspace{0.5cm}
\label{tab_6}
\renewcommand\arraystretch{1.5}
\scalebox{0.72}{
\begin{tabular}{|c|c|c|c|c|c|c|c|c|}
\hline
\multirow{3}{*}{$n$} &
\multicolumn{4}{c|}{$a_0$=0} &
\multicolumn{4}{c|}{$a_0$=1} \\
\cline{2-9} 
&\multicolumn{2}{c|}{The result of $u_\delta^{1}$} &
\multicolumn{2}{c|}{The result of $u_\delta^{2}$} &
\multicolumn{2}{c|}{The result of $u_\delta^{1}$} &
\multicolumn{2}{c|}{The result of $u_\delta^{2}$} \\
\cline{2-9}  
 & $\|u-u_n\|_{H^1}$ & $n^{-2/3}$ & $\|u-u_n\|_{H^1}$ & $n^{-2/3}$& $\|u-u_n\|_{H^1}$ & $n^{-2/3}$ & $\|u-u_n\|_{H^1}$ & $n^{-2/3}$  \\ \hline
 16 & 1.955e+00 & 0.00 & 1.955e+00 & 0.00 & 1.955e+00 & 0.00 &1.955e+00 & 0.00 \\ \hline 
32 & 1.950e+00 & 0.00 & 1.945e+00 & 0.01 & 1.950e+00 & 0.00 &1.943e+00 & 0.01 \\ \hline 
64& 1.616e+00 & 0.27 & 1.729e+00 & 0.17 & 1.607e+00 & 0.28 & 1.601e+00 & 0.28 \\ \hline
128  & 5.947e-01 & 1.44 & 5.569e-01 & 1.63 & 5.850e-01 & 1.46 & 5.326e-01 & 1.59 \\ \hline 
256 & 1.878e-01 & 1.66 & 1.821e-01 & 1.61 & 1.878e-01 & 1.64 & 1.842e-01 & 1.53 \\ \hline
\end{tabular}}
\end{table} 
\subsubsection{The case $k = 2$}
In this section, we test the case when $k = 2$ is used for the activation ${\rm ReLu}^k$ function. Note that the optimal rate is given as 
\begin{equation}
\|u - u_n\|_1 \approx n^{-1}.
\end{equation}
Thus, if we choose $\delta = n^{-\frac{2}{3}}$, then we expect that the optimal convergence rate can be achieved.

\begin{table}[H]
\centering
\caption{Numerical results of OGA for 3D elliptic equation with Dirichlet boundary condition and $\delta =256^{-\frac{2}{3}}$. We show two different results for $a_0=0$ and $a_0=1$. The predicted $H^1$ norm error is about $2\delta^{\frac{3}{2}} = $ 7.8e-03.} 
\vspace{0.5cm}
\label{tab_7}
\renewcommand\arraystretch{1.5}
\scalebox{0.72}{
\begin{tabular}{|c|c|c|c|c|c|c|c|c|}
\hline
\multirow{3}{*}{$n$} &
\multicolumn{4}{c|}{$a_0$=0} &
\multicolumn{4}{c|}{$a_0$=1} \\
\cline{2-9} 
&\multicolumn{2}{c|}{The result of $u_\delta^{1}$} &
\multicolumn{2}{c|}{The result of $u_\delta^{2}$} &
\multicolumn{2}{c|}{The result of $u_\delta^{1}$} &
\multicolumn{2}{c|}{The result of $u_\delta^{2}$} \\
\cline{2-9}  
 & $\|u-u_n\|_{H^1}$ &$n^{-1}$ & $\|u-u_n\|_{H^1}$ & $n^{-1}$& $\|u-u_n\|_{H^1}$ & $n^{-1}$ & $\|u-u_n\|_{H^1}$ & $n^{-1}$  \\ \hline
 16  & 1.955e+00 & 0.00  & 1.955e+00 & 0.00 & 1.955e+00 & 0.00& 1.956e+00 & 0.00 \\ \hline 
32 & 1.829e+00 & 0.10 & 1.864e+00 & 0.07 & 1.801e+00 & 0.12 & 1.868e+00 & 0.07 \\ \hline 
64  & 3.152e-01 & 2.54 &  4.664e-01 & 2.00 & 2.386e-01 & 2.92 & 4.210e-01 & 2.15 \\ \hline 
128  & 5.044e-02 & 2.64 &  5.773e-02 & 3.01 & 3.694e-02 & 2.69  & 5.499e-02 & 2.94 \\ \hline 
256  & 5.377e-03 & 3.23 & 6.410e-03 & 3.17 & 6.798e-03 & 2.44 &  4.598e-03 & 3.58 \\ \hline 
\end{tabular}}
\end{table} 

\section{Conclusions}
In this paper, we proposed an algorithm that can lead to the optimal error estimate of the shallow neural network approximation to the elliptic PDEs with the Dirichlet boundary condition. The newly proposed algorithm is shown to lead to the optimal error bound in $H^1$ norm with relatively small penalty parameter. The main technical tool is to use the duality argument to get an error estimate in $L^2$ norm and the universal approximation property. A sample numerical tests have been constructed to demonstrate the effectivity and optimal accuracy of the proposed method for 1D, 2D and 3D. 

\section*{Acknowledgments}
Authors would like to express their sincere appreciation to Professor Jinchao Xu for a number of his stimulating lectures on machine learning algorithms and guidance. 
\bibliographystyle{siamplain}
\bibliography{references}
\end{document}